\documentclass[psfig]{article}
\usepackage{amssymb, amsmath, amsthm}
\usepackage{MnSymbol}
\usepackage{esint}
\usepackage{wasysym}
\usepackage{graphicx}
\usepackage{color}
\usepackage{epsfig}
\usepackage{cite}

\theoremstyle{plain}
\newtheorem{thm}{Theorem}

\theoremstyle{definition}

\theoremstyle{remark}

\newcommand{\norm}[1]{\left\Vert#1\right\Vert}

\newcommand{\Real}{\mathbb R}

\begin{document}
\bibliographystyle{plain}

\title{Existence of Weak Solutions for the Incompressible Euler Equations}
\author{Emil Wiedemann}
\date{}
\maketitle
\begin{abstract}
Using a recent result of C. De Lellis and L. Sz\'{e}kelyhidi Jr.
(\cite{euler2}) we show that, in the case of periodic boundary
conditions and for arbitrary dimension $d\geq 2$, there exist
infinitely many global weak solutions to the incompressible Euler
equations with initial data $v_0$, where $v_0$ may be any solenoidal
$L^2$-vectorfield. In addition, the energy of these solutions is
bounded in time.
\end{abstract}

\section{Introduction}
Let $Q=[0,2\pi]^d$, $d\geq2$, and $L^2_{per}(Q)$ be the space of
$Q$-periodic functions in $L^2_{loc}(\Real^d;\Real^d)$, i.e.
$u(x+2\pi l)=u(x)$ for a.e. $x\in\Real^d$ and every
$l\in\mathbb{Z}^d$. Then, as usual when dealing with periodic
boundary conditions for fluid equations (cf. for instance
\cite{const}), we define the space
\begin{equation}
\begin{aligned}
H^m_{per}(Q)=&\{v\in L^2_{per}(Q):
\sum_{k\in\mathbb{Z}^d}|k|^{2m}|\hat{v}(k)|^2<\infty, \text{
}\hat{v}(k)\cdot k=0\text{   for every   }k\in\mathbb{Z}^d,\\
&\text{ and   }\hat{v}(0)=0 \},
\end{aligned}
\end{equation}
where $\hat{v}:\mathbb{Z}^d\rightarrow\mathbb{C}^d$ denotes the
Fourier transform of $v$. We shall write $H(Q)$ instead of
$H^0_{per}(Q)$ and $H_w(Q)$ for the space $H(Q)$ equipped with the
weak $L^2$ topology.

Recall the incompressible Euler equations
\begin{equation}
\begin{aligned}
\partial_tv+\operatorname{div}(v\otimes v)+\nabla p&=0\\
\operatorname{div}v&=0,
\end{aligned}
\end{equation}
where $v\otimes v$ is the matrix with entries $v_iv_j$ and the
divergence is taken row-wise. A vectorfield $v\in
L^{\infty}\left((0,\infty);H(Q)\right)$ is called a \emph{weak
solution} of these equations with $Q$-periodic boundary conditions
and initial data $v_0\in H(Q)$ if
\begin{equation}
\int_0^{\infty}\int_{Q}(v\cdot\partial_t\phi+v\otimes
v:\nabla\phi)dxdt+\int_{Q}v_0(x)\phi(x,0)dx=0
\end{equation}
for every $Q$-periodic divergence-free $\phi\in
C_c^{\infty}\left(\Real^d\times[0,\infty);\Real^d\right)$.

Unlike in the case of Navier-Stokes equations, for which the
existence of global weak solutions has been known since the work
\cite{leray} of J. Leray, the existence problem for weak solutions
of Euler has remained open so far. In this paper we show that the
existence of weak solutions is a consequence of C. De Lellis' and L.
Sz\'{e}kelyhidi's work \cite{euler2}. More precisely, we have
\begin{thm}\label{existence}
Let $v_0\in H(Q)$. Then there exists a weak solution $v\in
C([0,\infty);H_w(Q))$ (in fact, infinitely many) of the Euler
equations with $v(0)=v_0$. Moreover, the kinetic energy
\begin{equation}
E(t):=\frac{1}{2}\int_Q|v(x,t)|^2dx
\end{equation}
is bounded and satisfies $E(t)\rightarrow0$ as $t\rightarrow\infty$.
\end{thm}
Note that the condition $\hat{v}(0)=0$ in the definition of $H(Q)$,
i.e. $\int_Qvdx=0$, is no actual constraint due to Galilean
invariance of the Euler equations.

Our proof of this theorem is very simple: Owing to \cite{euler2}, it
suffices to construct a suitable so-called \emph{subsolution} with
the desired initial data; we obtain such a subsolution by solving
the Cauchy problem for the fractional heat equation
\begin{equation}
\begin{aligned}
\partial_tv+(-\Delta)^{1/2}v&=0\\
v(\cdot,0)&=v_0,
\end{aligned}
\end{equation}
which is not difficult since, owing to periodicity, we can work in
Fourier space.

Although our solutions have bounded energy, they do not satisfy any
form of the energy inequality. Indeed, they exhibit an increase in
energy at least at time $t=0$, and this increase will be
discontinuous (this follows from $e(v_0,u_0)>\frac{1}{2}|v_0|^2$ in
the proof below). If one requires, in contrast, that the energy be
bounded at all times by the initial energy, then existence of such
weak solutions is not known for arbitrary initial data (but only for
an $L^2$-dense subset of initial data, see \cite{generalyoung}). In
fact, it is impossible to deduce from Theorem \ref{convexint} below
such an existence theorem, since for smooth initial data the
existence of infinitely many weak solutions would contradict
well-known local existence results and weak-strong uniqueness, see
Subsection 2.3 of \cite{euler2}.

\section{Preliminaries}
Before we prove the result of this paper, we recall some notions
from \cite{euler2}. Let $\mathcal{S}_0^d$ denote the space of
symmetric trace-free $d\times d$-matrices. Then the
\emph{generalised energy}
$e:\Real^d\times\mathcal{S}_0^d\rightarrow\Real$ is defined by
\begin{equation}
e(v,u)=\frac{n}{2}\lambda_{max}(v\otimes v-u),
\end{equation}
where $\lambda_{max}$ denotes the largest eigenvalue. $e$ is known
to be non-negative and convex, and $\frac{1}{2}|v|^2\leq e(v,u)$ for
all $v$ and $u$ with equality if and only if $u=v\otimes
v-\frac{|v|^2}{d}I_d$ ($I_d$ being the $d\times d$ unit matrix). The
following is shown in \cite{euler2}:
\begin{thm}\label{convexint}
Let $\bar{e}\in C\left(\Real^d\times(0,\infty)\right)\cap
C\left([0,\infty);L^1_{loc}(\Real^d)\right)$ be $Q$-periodic in the
space variable and $(\bar{v},\bar{u},\bar{q})$ be a smooth,
$Q$-periodic (in space) solution of
\begin{equation}\label{linear}
\begin{aligned}
\partial_t\bar{v}+\operatorname{div}\bar{u}+\nabla\bar{q}&=0\\
\operatorname{div}\bar{v}&=0
\end{aligned}
\end{equation}
in $\Real^d\times(0,\infty)$ such that
\begin{equation}
\bar{v}\in C([0,\infty);H_w(Q)),
\end{equation}
\begin{equation}
\bar{u}(x,t)\in\mathcal{S}_0^d
\end{equation}
for every $(x,t)\in Q\times(0,\infty)$, and
\begin{equation}
e\left(\bar{v}(x,t),\bar{u}(x,t)\right)<\bar{e}(x,t)
\end{equation}
for every $(x,t)\in Q\times(0,\infty)$.

Then there exist infinitely many weak solutions $v\in
C([0,\infty);H_w(Q))$ of the Euler equations with
$v(x,0)=\bar{v}(x,0)$ for a.e. $x\in Q$ and
\begin{equation}
\frac{1}{2}|v(x,t)|^2=\bar{e}(x,t)
\end{equation}
for every $t\in(0,\infty)$ and a.e. $x\in Q$.
\end{thm}

\section{Proof of Theorem \ref{existence}}
By Theorem \ref{convexint}, it suffices to find suitable
$(\bar{v},\bar{u},\bar{q})$ and $\bar{e}$.

Let us define $\bar{v}$ and $\bar{u}$ by their Fourier transforms as
follows:
\begin{equation}\label{defv}
\hat{\bar{v}}(k,t)=e^{-|k|t}\hat{v}_0(k),
\end{equation}
\begin{equation}\label{defu}
\hat{\bar{u}}_{ij}(k,t)=-i\left(\frac{k_j}{|k|}\hat{\bar{v}}_i(k,t)+\frac{k_i}{|k|}\hat{\bar{v}}_j(k,t)\right)
\end{equation}
for every $k\neq0$, and $\hat{\bar{u}}(0,t)=0$. Note that
$\bar{u}_{ij}$ thus defined equals
$-\mathcal{R}_j\bar{v}_i-\mathcal{R}_i\bar{v}_j$, where
$\mathcal{R}$ denotes the Riesz transform. Clearly, for $t>0$,
$\bar{v}$ and $\bar{u}$ are smooth. Moreover, $\bar{u}$ is symmetric
and trace-free. Indeed, the latter can be seen by observing
\begin{equation}
\sum_{i=1}^d\left(\frac{k_i}{|k|}\hat{\bar{v}}_i(k,t)+\frac{k_i}{|k|}\hat{\bar{v}}_i(k,t)\right)=\frac{2}{|k|}e^{-|k|t}k\cdot\hat{v}_0(k)=0
\end{equation}
for all $k\neq0$ (for $k=0$ this is obvious).

Next, we can write equations (\ref{linear}) in Fourier space as
\begin{equation}\label{linearfourier}
\begin{aligned}
\partial_t\hat{\bar{v}}_i+i\sum_{j=1}^dk_j\hat{\bar{u}}_{ij}+ik_i\hat{\bar{q}}&=0\\
k\cdot\hat{\bar{v}}&=0
\end{aligned}
\end{equation}
for $k\in\mathbb{Z}^d$, $i=1,\ldots,d$. It is easy to check that
$(\hat{\bar{v}},\hat{\bar{u}},0 )$ as defined by (\ref{defv}) and
(\ref{defu}) solves (\ref{linearfourier}) and hence
$(\bar{v},\bar{u},0)$ satisfies (\ref{linear}).

Concerning the energy, we have the pointwise estimate
$e(\bar{v},\bar{u})\leq C(|\bar{v}|^2+|\bar{u}|)$, and because of
\begin{equation}
\int_Q|\bar{v}|^2dx=\sum_{k\in\mathbb{Z}^d}|\hat{\bar{v}}|^2=\sum_{k\in\mathbb{Z}^d}e^{-|k|t}|\hat{v}_0|^2\leq\norm{v_0}^2_{L^2(Q)}
\end{equation}
and, similarly,
\begin{equation}
\int_Q|u|dx\leq C\int_Q|u|^2dx\leq C\norm{v_0}^2_{L^2(Q)},
\end{equation}
we conclude that
$\sup_{t>0}\norm{e(\bar{v}(x,t),\bar{u}(x,t))}_{L^1(Q)}<\infty$.
Moreover, from the same calculation and the dominated convergence
theorem we deduce
\begin{equation}
\norm{e(\bar{v}(x,t),\bar{u}(x,t))}_{L^1(Q)}\rightarrow0
\end{equation}
as $t\rightarrow\infty$ as well as
\begin{equation}
\bar{v}(t)\rightarrow v_0
\end{equation}
strongly in $L^2(Q)$ and
\begin{equation}
\bar{u}(t)\rightarrow
u_0:=-(\mathcal{R}_j(v_0)_i+\mathcal{R}_i(v_0)_j)_{ij}
\end{equation}
strongly in $L^1(Q)$. We claim that then
\begin{equation}
e(\bar{v},\bar{u})\in C\left([0,\infty);L^1(Q)\right).
\end{equation}
The only issue is continuity at $t=0$. First, one can easily check
that the map
\begin{equation}
(v,u)\mapsto e\left(\frac{v}{\sqrt{|v|}},u\right)
\end{equation}
is Lipschitz continuous with Lipschitz constant, say, $L$; thus,
using the inequality $||a|a-|b|b|\leq(|a|+|b|)|a-b|$, we have
\begin{equation}
\begin{aligned}
\int_Q|e(\bar{v},\bar{u})-e(v_0,u_0)|&\leq
L\int_Q\left(||\bar{v}|\bar{v}-|v_0|v_0|+|\bar{u}-u_0|\right)dx\\
&\leq
2L\sup_{t\geq0}\norm{\bar{v}(t)}_{L^2}\norm{\bar{v}(t)-v_0}_{L^2}+L\norm{\bar{u}-u_0}_{L^1}\rightarrow0
\end{aligned}
\end{equation}
as $t\rightarrow0$. This proves the claim.

Therefore, $\bar{e}$ defined by
\begin{equation}
\bar{e}(x,t):=e(\bar{v}(x,t),\bar{u}(x,t))+\min\{t,\frac{1}{t}\}
\end{equation}
satisfies the requirements of Theorem \ref{convexint} and, in
addition, $\int_Q\bar{e}dx\rightarrow0$ as $t\rightarrow\infty$.
Theorem \ref{convexint} then yields the desired weak solutions of
Euler. \qed

\bibliography{EulerExistence}

\end{document}